\documentstyle[12pt,verbatim]{amsart}

\numberwithin{equation}{section}

\renewcommand{\thesection}{\arabic{section}}
\renewcommand{\theequation}{\thesection.\arabic{equation}}

\newtheorem{theorem}[equation]{Theorem}
\newtheorem{lemma}[equation]{Lemma}
\newtheorem{corollary}[equation]{Corollary}
\newtheorem{proposition}[equation]{Proposition}

\theoremstyle{definition}

\theoremstyle{remark}

\newcounter{itemno}
\newenvironment{num}{\mbox{}\begin{list}{(\theequation.\theitemno)}{\settowidth{\labelwidth}{(M.1M.M)}\usecounter{itemno}}}{\end{list}}

\newcommand{\NS}{\operatorname{NS}}
\newcommand{\MW}{\operatorname{MW}}
\newcommand{\rank}{\operatorname{rank}}
\newcommand{\Hom}{\operatorname{Hom}}

\def\B{ {\bold{B}} }
\def\M{ {\bold{M}} }
\def\P{ {\bold{P}} }
\def\Re{ {\bold{Re}} }

\def\1{ {\bold{1}} }

\begin{document}
\title{Torsion Sections of Semistable Elliptic Surfaces}
\author{Rick Miranda}
\address{Dept. of Mathematics\\Colorado State University\\Ft. Collins, CO
80523}
\email{miranda@@riemann.math.colostate.edu}
\thanks{Research supported in part by the NSF under grant DMS-9104058\\ \today}

\keywords{elliptic surfaces, Mordell-Weil, sections, torsion}
\subjclass{14J27 secondary 14H52}
\maketitle


\setcounter{section}{-1}

\section{Introduction}

In this article we would like to address the following situation.
Let $f:X \to C$ be an elliptic surface with section $S_0$.
We further assume that $X$ is relatively minimal and smooth,
and that all singular fibers are semistable,
that is, are cycles of ${\Bbb P}^1$'s
(type ``$I_m$'' in Kodaira's notation, see \cite{kodaira}).
Given a section $S$ of $f$,
it will meet one and only one component in every singular fiber.
The fundamental aim of this paper is
to try to say something specific
about which components are being hit,
in the case that $S$ is a torsion section
(i.e., of finite order in the Mordell-Weil group of sections of $f$).

Most work on the Mordell-Weil group for elliptic surfaces
focuses on properties of the group itself,
(for example its rank, etc.)
and not on properties of the elements of the group.
See for example \cite{shioda} and \cite{cox-zucker} for general properties,
and papers such as \cite{schwartz} and \cite{stiller}
for computations in specific cases.
Although the methods used in this article are elementary,
they are suprisingly powerful,
enabling us to obtain rather detailed information
about the possibilities for component intersection
for a torsion section.

To be specific, suppose that the $j^{th}$ singular fiber
is a cycle of length $m_j$, and that the components
are numbered cyclically around the cycle from $0$ to $m_j-1$.
Suppose that in this $j^{th}$ fiber,
a section $S$ of order exactly $n$ hits the $k_j^{th}$ component.
The main restrictions on these ``component'' numbers $k_j$ are the following:

\medskip
\noindent
{\bf The Quadratic Relation} (Proposition \ref{formula1}):
\newline
\[
\sum_j k_j(1 - k_j/m_j) = 2 \chi({\cal O}_X).
\]
\noindent
{\bf The Component Number Sum Relation}
(Corollaries \ref{formula2} and \ref{formula4}):
\newline
If the component numbers satisfy $k_j \leq m_j/2$,
(which can always be arranged by the suitable ``reorienting'' of the cycle),
then
\[
\sum_{j=1}^s k_j(S) = \begin{cases}
4\chi({\cal O}_X) & \text{ if } n=2 \\
3\chi({\cal O}_X) & \text{ if } n \geq 3.
\end{cases}
\]

These relations hold for the general torsion section,
but actually we have more to say in the case
of a torsion section $S$ of prime order $p$.
In this case we are able to calculate exactly the distribution
of the ``distances'' between the components hit by $S$ and the components
hit by the zero-section $S_0$, in the following sense.
If a fiber has $m$ components, indexed from $0$ to $m-1$ cyclically,
then $S$, having order $p$, meets either the component indexed by $0$
(which we take to be the component meeting $S_0$)
or one of the components indexed by $im/p$ for some $i$.
In other words, a section of order $p$ must meet either $1/p$ of the
way around the cycle, or $2/p$ of the way around, etc.
Denote by $M_{i,p}(S)$ the fraction of the total Euler number $e=\sum_j m_j$
contributed by those fibers in which $S$ meets $i/p$ of the way around
the cycle, in either direction.
Our main result is a computation of these numbers,
which turn out to be independent of $S$, and surprisingly,
independent of $i$ (as long as $i\neq 0$):

\medskip
\noindent
{\bf The Distribution Numbers} (Theorem \ref{Minvalues}):
\newline
Let $S$ be a torsion section of odd prime order $p$.
Then
\[
M_{i,p}(S) = 2p/(p^2-1) \; \text{ if }\; 1 \leq i \leq (p-1)/2,
\]
and
\[
M_{0,p}(S) = 1/(p+1).
\]

Thus we obtain an ``equidistribution'' of these distances:
the total sum of those $m_j$'s contributed from fibers of type $I_{m_j}$
in which $S$ meets $i/p$ of the way around the cycle is constant,
independent of $i$.

Thus a torsion section of prime order,
as it travels across the elliptic surface,
is forced to be completely ``fair'' in its choit may hit, in the above sense.

The above results generalize various statements
concerning the component numbers $\{k_j\}$
for a torsion section
which appeared in \cite{miranda-persson1},
\cite{miranda-persson2}, \cite{miranda-persson3},
and \cite{miranda}.

The results in this paper could be proved
by appealing to the appropriate elliptic modular surface
(for those primes where such exists)
and its universal properties.
The equidistribution property is invariant under base change,
so if it is true for the (universal) modular surface,
it will be true in general.
This is a computation involving the modular group,
which could be made once and for all.
However it is my intention to show that these properties are somehow
more elementary than the theory of modular surfaces,
and follow from straightforward considerations involving only
basic facts concerning the Mordell-Weil group
and a little intersection theory of surfaces.

In a further paper we hope to discuss the case of torsion
of non-prime order more fully.

I would like to thank Igor Dolgachev, Jeanne Duflot, Ulf Persson,
Miles Reid, and Peter Stiller for some useful conversations.
I would also like to thank the organizers of the Conference
in L'Aquila for their kind invitation to allow me to speak on this subject.
I would especially like to thank Prof. E. Laura Livorni
for her wonderful efforts to make the conference such a success,
and Prof. Ciro Ciliberto for his excellent hospitality in Rome
while this work was completed.

\section{Basic Facts}

We work over the field of complex numbers.
Let $C$ be a smooth curve, and
let $f:X \to C$ be a smooth semistable elliptic surface over $C$.
We suppose that $f$ has a chosen section $S_0$
which we treat as the zero for the group law on the sections of $f$.
Let $g$ be the genus of $C$,
and let $\chi = \chi({\cal O}_X)$
be the holomorphic Euler characteristic of $X$.
If $e$ is the topological Euler number of $X$, then $e = 12 \chi$.

We suppose that the map $f$ has $s$ singular fibers $F_1,\dots,F_s$,
each semistable, i.e., each a cycle of rational curves
(of type ``$I_m$'' in Kodaira's notation \cite{kodaira}).
Indeed, let us say that the fiber $F_j$ is of type $I_{m_j}$.
Therefore
\[
e = 12 \chi = \sum_{j=1}^s m_j.
\]

Choose an ``orientation'' of each fiber $F_j$ and write the $m_j$ components
of $F_j$ as
\[
C^{(j)}_0,C^{(j)}_1,\dots,C^{(j)}_{m_j-1},
\]
where the zero section $S_0$ meets only $C^{(j)}_0$
and for each $k$, $C^{(j)}_k$ meets only $C^{(j)}_{k\pm 1} \mod m_j$.
If $m_j = 1$, then $F_j = C^{(j)}_0$ is a nodal rational curve,
of self-intersection $0$.
If $m_j \geq 2$, then each $C^{(j)}_k$ is a smooth rational curve
with self-intersection $-2$.

If we denote by $\NS(X)$ the Neron-Severi group of the elliptic surface $X$,
we can consider the sublattice $L \subseteq \NS(X)$ generated by the
zero-section $S_0$ and the components $C^{(j)}_k$ of the singular fibers.
The class of the fiber $F$ is of course in this sublattice $L$;
moreover, the only relations among these classes are that
\begin{equation}
\label{Frelation}
F \equiv \sum_{k=0}^{m_j-1} C^{(j)}_k \text{ for each } j = 1,\dots,s.
\end{equation}
If we let $U$ denote the sublattice of $L$ spanned by $F$ and $S_0$,
then $U$ is a unimodular lattice of rank $2$,
isomorphic to a hyperbolic plane, and so splits off $L$.
Its orthogonal complement $R$ is freely generated by the components
$C^{(j)}_k$ for $j = 1 \dots s$ and $k \neq 0$.
$R$ is therefore a direct sum of $s$ root lattices of type $A$,
with the $j^{th}$ summand isomorphic to $A_{m_j-1}$,
generated by $C^{(j)}_k$ for $k \neq 0$.

The Shioda-Tate formula for the Mordell-Weil group $\MW(X)$
of sections of $X$
is derived from the exact sequence (see \cite{shioda})
\[
0 \to L \to \NS(X) \to \MW(X) \to 0
\]
where the first map is the inclusion of the sublattice $L$ into $\NS(X)$.
The second map is the fiber-by-fiber summation map,
sending a divisor class $D \in \NS(X)$
to the closure of the sum of the points of $D$ on the generic fiber of $X$.
We therefore obtain information about both the rank and the torsion
of the Mordell-Weil group $\MW(X)$.  If we denote by $\rho$ the Picard number
of $X$, which is the rank of the Neron-Severi group, we see that since
$L$ has rank equal to $2 + \sum_{j=1}^s (m_j-1)$,
\[
\rank\MW(X) = \rho - 2 - \sum_{j=1}^s (m_j-1) = s + \rho - 2 - e.
\]
The torsion $\MW_{tor}(X)$ in the Mordell-Weil group
corresponds to those classes in $\NS(X)$
for which some multiple lies in the sublattice $L$.
These classes form an intermediate lattice $L^{sat}$
containing $L$ as a sublattice of finite index,
and we see that
\[
\MW_{tor}(X) \cong L^{sat}/L.
\]

For any lattice $N$,
denote by $N^\#$ the dual lattice $\Hom_{\Bbb Z}(N,\Bbb Z)$.
The inclusion of $L$ into $L^{sat}$ gives the sequence of inclusions
\[
L \subseteq L^{sat} \subseteq {(L^{sat})}^\# \subseteq L^\#
\]
which shows that the quotient $L^{sat}/L$ is isomorphic to
a subgroup of the discriminant-form group $G_L = L^\#/L$.
Note that $L^\# = U^\#\oplus R^\#$ since $L = U \oplus R$
and since $U^\# = U$ ($U$ is unimodular),
we have $G_L \cong R^\#/R = \bigoplus_{j=1}^s A_{m_j-1}^\#/A_{m_j-1}$.
A computation shows that for the root lattices of type $A$,
$A_m^\#/A_m \cong \Bbb Z/m\Bbb Z$;
therefore we have that the torsion part $\MW_{tor}(X)$ of the Mordell-Weil
group of sections of $X$ is isomorphic to a subgroup of
$\bigoplus_{j=1}^s \Bbb Z/m_j\Bbb Z$.
(See for example \cite{miranda-persson1} or \cite{miranda}.)

A choice of orientation for each singular fiber $F_j$
gives an identification of the cyclic group of components
$C^{(j)}_k$ for $k =0,\dots,m_j-1$
with $\Bbb Z/m_j\Bbb Z$.
For any section $S$ of $f$,
and any singular fiber $F_j$,
denote by $k_j(S)$ the index of the component of $F_j$ which $S$ meets.
Thus
\begin{equation}
\label{Sdot}
S \cdot C^{(j)}_l = \delta_{l k_j(S)}
\end{equation}
where $\delta$ above is the Kronecker delta.
Note that with our notation above, $k_j(S_0) = 0$ for every $j$.
The numbers $\{k_j(S)\}$ will be called the {\em component numbers}
of the section $S$.
The upshot of the remarks above is that, for torsion sections,
this assignment of integers $\{k_j(S)\}$ to $S$ is $1$-$1$.
Thus a torsion section is determined by its
component numbers (but not all sets of components numbers can occur).

It is the purpose of this article to study
the properties of the component numbers of torsion sections
for semistable elliptic surfaces.
Although it is clear from the above discussion
that the component numbers $k_j(S)$ are well-defined modulo $m_j$,
it is more useful for our purposes
to take them to be integers in the range $0,\dots, m_j-1$.

Two facts are necessary for what follows;
they can be found in \cite{miranda} as well as the other basic references.
\begin{lemma}
\label{2facts}
\begin{num}
\item If $S_1$ and $S_2$ are two different sections in $\MW(X)$
with $S_1-S_2$ torsion, then $S_1\cdot S_2 = 0$, i.e., they are disjoint.
\item If $S$ is any section in $\MW(X)$, then $S\cdot S = -\chi$.
\end{num}
\end{lemma}

\section{The Divisor Class of a Torsion Section}

Let $S$ be a torsion section of the semistable elliptic surface
$f:X \to C$.
In this section we wish to write down the class of $S$ in $\NS(X)$.

Since the class of $S$ lies in $L^{sat}$,
we have that $S$ is a $\Bbb Q$-linear combination
of the zero section $S_0$, and the fiber components $C^{(j)}_k$.
Because (\ref{Frelation}) is the only relation among these classes,
the classes $S_0$, $F$, and $C^{(j)}_k$ for $j = 1,\dots,s$
and $k \neq 0$ form a basis for the module spanned by all these classes.

For fixed $j=1,\dots,s$ and $k = 1,\dots,m_j-1$, set
\begin{eqnarray*}
D^{(j)}_k &=& (m_j - k) \sum_{i=1}^k i C^{(j)}_i +
                    k \sum_{i=k+1}^{m_j-1} (m_j - i) C^{(j)}_i \\
&=& (m_j-k)[ C^{(j)}_1 + 2C^{(j)}_2 + \dots + kC^{(j)}_k ] + \\
& & k [ (m_j-k-1)C^{(j)}_{k+1} + \dots + 2C^{(j)}_{m_j-2} + C^{(j)}_{m_j-1} ].
\end{eqnarray*}
We set $D^{(j)}_0 = 0$.

Note that
\begin{equation}
\label{S0dot}
D^{(j)}_k \cdot S_0 = 0
\end{equation}
for every $k$.

\begin{lemma}
\label{Dlemma}
Fix $j = 1,\dots,s$, and indices $k$ and $l$ with $0 \leq k,l \leq m_j-1$.
Then
\[
D^{(j)}_k \cdot C^{(j)}_l =
\begin{cases}
m_j & \text{ if } l = 0 \text{ and } k \neq 0 \\
-m_j & \text{ if } l = k \text{ and } k \neq 0 \\
0 & \text{ otherwise.}
\end{cases}
\]
\end{lemma}

\begin{pf}
This is a straightforward check;
let us first notice that if $k = 0$, then $D^{(j)}_k = 0$,
so certainly the intersection product is $0$.
Assume then that $k \neq 0$.
If $l=0$, then since $C^{(j)}_0$ meets only $C^{(j)}_1$ and $C^{(j)}_{m_j-1}$,
each once, we have
\[
D^{(j)}_k \cdot C^{(j)}_0 = (m_j-k)[1] + (k)[1] = m_j,
\]
which proves the first statement.  If $l = k$, then
$C^{(j)}_k$ meets only $C^{(j)}_{k\pm 1}$ (each once) and itself (-2 times),
so that
\[
D^{(j)}_k \cdot C^{(j)}_k = (m_j-k)[(k-1)(1) + (k)(-2)] + (k)[(m_j-k-1)(1)] =
-m_j,
\]
proving this case.
We leave to the reader to check that all other intersection products are $0$,
as stated.
\end{pf}

It is useful to note that the above lemma can be re-expressed
using the Kronecker delta function as
\begin{equation}
\label{Ddot}
D^{(j_1)}_k \cdot C^{(j_2)}_l =
m_j \delta_{j_1j_2} (1-\delta_{k0})(\delta_{l0} - \delta_{lk}).
\end{equation}

\begin{theorem}
With the above notation, if $S$ is a torsion section of order $n$,
then
\[
S \equiv S_0 + (S - S_0 \cdot S_0) F
- \sum_{j=1}^s {1 \over m_j} D^{(j)}_{k_j(S)}.
\]
\end{theorem}

\begin{pf}
Since the intersection form on $L$ is nondegenerate,
we need only check that both sides of the equation intersect the generators
$S_0$, $F$, and the $C^{(j)}_l$ in the same number.
Let us begin with $S_0$. Using (\ref{S0dot}) and the equations
\[
S \cdot F = S_0 \cdot F = 1,
\]
we have that the right hand side intersects $S_0$ to $S\cdot S_0$,
agreeing with the left hand side.
Similarly, $F$ meets only the $S_0$ term on the right hand side,
so the intersections with $F$ are also equal.

Now fix indices $j = 1,\dots,s$ and $l = 0,\dots, m_j-1$,
and let us check the intersection with $C^{(j)}_l$;
Lemma \ref{Dlemma} is the critical part of the computation.
Intersecting with the right hand side gives
$S_0 \cdot C^{(j)}_l - {1 \over m_j} D^{(j)}_{k_j(S)} \cdot C^{(j)}_l$,
which reduces to
$\delta_{l 0} - (1-\delta_{k_j(S) 0})(\delta_{l 0} - \delta_{l k_j(S)})$
using (\ref{Ddot}).
This is after some simplification equal to $\delta_{l k_j(S)}$,
which is also $S \cdot C^{(j)}_l$ by (\ref{Sdot}).
\end{pf}

The reader familiar with \cite{cox-zucker} will recognize
the above expression as the ``correction'' term
in their pairing on the Mordell-Weil group.
They were more interested in questions of rank in that paper,
and less so in the torsion sections.

The formula for the linear equivalence class of $S$ given above
can be simplified somewhat in the case that $S$ is not the zero section $S_0$.
Then $(S\cdot S_0)=0$ and $(S_0\cdot S_0) = -\chi$ by Lemma \ref{2facts},
so we have

\begin{equation}
\label{nonzeroS}
S \equiv S_0 + \chi F - \sum_{j=1}^s {1 \over m_j} D^{(j)}_{k_j(S)}.
\end{equation}

\section{The Quadratic Relation for the Component Numbers}

Let us now take the above formula for the torsion section $S$
and intersect it with $S$ itself.
If $S = S_0$, then of course all $k_j(S) = 0$,
both sides of the equation are $S_0$,
so we recover no information.
However if $S \neq S_0$,
the formula (\ref{nonzeroS}) is nontrivial.
Dotting $S$ with the right-hand side gives
\[
\chi - \sum_{j=1}^s {1 \over m_j} D^{(j)}_{k_j(S)} \cdot S.
\]
Since $S$ meets only the curve $C^{(j)}_{k_j(S)}$ in the $j^{th}$
singular fiber, and it meets it exactly once, we have
\[
D^{(j)}_{k_j(S)} \cdot S = (m_j - k_j(S))k_j(S).
\]
Therefore the above reduces to
\[
\chi - \sum_{j=1}^s k_j(S) (1-{k_j(S) \over m_j}).
\]
Finally, using Lemma \ref{2facts}.2,
dotting the left-hand side with $S$ gives $-\chi$.
Hence we obtain the following formula,
which we call the {\em quadratic relation} for the component numbers.

\begin{proposition}
\label{formula1}
Let $S$ be a torsion section of $f$, not equal to the zero section $S_0$.
Then
\[
\sum_{j=1}^s k_j(S) (1-{k_j(S) \over m_j}) = 2\chi({\cal O}_X).
\]
\end{proposition}

Note that the quadratic relation is independent of
the choice of orientation of each singular fiber, as it should be.
(If one reverses the orientation of $F_j$,
then $k_j$ is replaced by $m_j-k_j$ if $k_j \neq 0$,
so that the two terms of each summand are simply switched.)

Miles Reid has pointed out to me that the quadratic relation given
above can be derived directly via the use of the Riemann-Roch
theorem for Weil divisors, as expounded in \cite[section 9]{reid}.
One should contract all components of fibers except $C_0^{(j)}$,
obtaining a surface with only $A_{m_j-1}$ singularities,
and consider on that surface the Weil divisor $S$.
The Riemann-Roch formula for $S$ has as corrections to the usual
Riemann-Roch the terms $k_j(S) (1-{k_j(S) \over m_j})$.
(See also \cite[section 2]{giraud}.)
This approach is cleaner and perhaps conceptually simpler;
however in the spirit of keeping things as elementary as possible,
I have decided to present this ``low road'' version.
The interested reader will have no problem making the computation
as Reid suggests and re-deriving the quadratic relation in this manner.

\section{The Component Number Sums}

In this section we will develop formulas for the component number sums
$\sum_j k_j(S)$ for a torsion section $S$.
Such formulas follow rather directly
from the quadratic relation for the component numbers.
Suppose first that $S$ has order $2$.
Then for each $j$ with $k_j(S) \neq 0$,
we must have $m_j$ even and $k_j(S) = m_j/2$.
Therefore the quadratic relation reduces to the following.

\begin{corollary}
\label{formula2}
Let $S$ be an order $2$ section of $f$.  Then
\[
\sum_{j=1}^s k_j(S) = 4\chi({\cal O}_X).
\]
\end{corollary}

Note that in the order $2$ case,
each $k_j$ is either $0$ or $m_j/2$,
and hence the component numbers for a torsion section of order $2$
are independent of the orientation of the fibers.
This is clearly not true for other values of the component numbers.
In general, reversing the orientation of the fiber $F_j$
changes the component number from $k_j$ to $m_j-k_j$.
Hence the analogue of Corollary (\ref{formula2})
for torsion sections of order at least $3$
must take this into account.

To this end define a function $d_m:\{0,\dots,m-1\} \to \{0,\dots,[m/2]\}$
by setting
\[
d_m(k) = \min\{k, m-k\}.
\]
Note that for any given torsion section $S$,
it is possible to choose the orientations of the fibers
so that the component numbers $k_j(S)$ are minimal,
that is, $k_j(S) = d_{m_j}(k_j(S))$ for each $j$.
If this condition holds,
we say that the section $S$ has {\em minimal component numbers}.

Next suppose that a torsion section $S$ has order $n \geq 3$.
Then both $S$ and $2S$
are nonzero torsion sections of $X$.
If for each index $j$, we choose the orientation of the components
so that $0 \leq k_j(S) \leq m_j/2$,
then we have the formulas
\begin{equation}
\label{kjfor2S}
k_j(2S) = \begin{cases}
0 \text{ if }k_j(S) = m_j/2, \text{ and } \\
2k_j(S) \text{ if }k_j(S) < m_j/2
\end{cases}
\end{equation}
for every $j$.
Therefore applying Proposition \ref{formula1} to $2S$,
after dividing by $2$ we obtain
\[
\sum\begin{Sb} j \text{ with }\\ k_j(S)<m_j/2 \end{Sb}
k_j(S) (1-2{k_j(S) \over m_j}) = \chi({\cal O}_X).
\]
Multiplying the formula of Proposition \ref{formula1} by $2$ gives
\[
\sum\begin{Sb} j \text{ with }\\ k_j(S)=m_j/2 \end{Sb} k_j(S)  +
\sum\begin{Sb} j \text{ with }\\ k_j(S)<m_j/2 \end{Sb}
k_j(S) (2-2{k_j(S) \over m_j}) = 4\chi({\cal O}_X),
\]
and subtracting the above two equations yields the following.

\begin{corollary}
\label{formula3}
Let $S$ be a section of order $n \geq 3$
with minimal component numbers $\{k_j(S)\}$.
Then
\[
\sum_{j=1}^s k_j(S) = 3\chi({\cal O}_X).
\]
\end{corollary}

In case the orientations of the fibers are not chosen so as
to give $S$ minimal component numbers, a similar formula holds,
expressed with the $d_m$ function.

\begin{corollary}
\label{formula4}
Let $S$ be a section of order $n \geq 3$.
Then
\[
\sum_{j=1}^s d_{m_j}(k_j(S)) = 3\chi({\cal O}_X).
\]
\end{corollary}

Using these component number sum relations,
one can re-express the quadratic relation
in a simpler form.

\begin{corollary}
\label{formula5}
Let $S$ be a section of order $n \geq 2$,
with minimal component numbers $\{k_j=k_j(S)\}$.
Then
\[
\sum_{j=1}^s k_j^2/m_j = \begin{cases}
\chi({\cal O}_X) & \text{ if } n \geq 3 \\
2\chi({\cal O}_X) & \text{ if } n = 2
\end{cases}
\]
\end{corollary}

\begin{pf}
The quadratic relation can be written as
\[
\sum_j k_j^2/m_j = \sum_j k_j - 2\chi({\cal O}_X),
\]
and so the result follows from Corollaries \ref{formula2} and \ref{formula4}
\end{pf}

These results generalize a certain divisibility statement concerning
the component numbers obtained in \cite[Proposition 4.6]{miranda-persson1}
(see also \cite[Lemma 3(a)]{miranda-persson2} and \cite[X.3.1]{miranda}).
Although the above cited works deal with elliptic K3 surfaces,
the statement holds in general.  It is:
\[
\sum_{j=1}^s {m_j-1 \over 2m_j} k_j^2 \; \text{ is an integer.}
\]
This was proved using lattice-theoretic arguments involving
the discriminant-form group $L^{\#}/L$.

We can recover this result as follows.
First note that in the above sum,
if any $k_j$ is replaced by $m_j-k_j$,
the sum changes by an integer;
hence it suffices for us to derive the divisibility result
using minimal component numbers.

Now the above sum in question can be rewritten as
\[
\sum_j k_j^2/2  -  \sum_j  k_j^2/2m_j.
\]
Note that $\sum_j k_j^2$ has the same parity as $\sum_j k_j$,
so that the first term $\sum_j k_j^2/2$ is equal to
$\sum_j k_j/2$ modulo ${\Bbb Z}$.
Now modulo ${\Bbb Z}$, we have

\begin{eqnarray*}
\sum_j k_j^2/2  -  \sum_j  k_j^2/2m_j &=& \sum_j k_j/2  -  \sum_j  k_j^2/2m_j
\\
&=& \begin{cases}
4\chi/2 - 2\chi/2  & \text{ if } n=2 \\
3\chi/2 - \chi/2 & \text{ if } n \geq 3
\end{cases} \\
&=& \chi
\end{eqnarray*}
which is an integer.

\section{The Distribution Numbers for a Torsion Section}

In this section we will apply the relations
developed in the previous sections
for the component numbers of a non-zero torsion section $S$
to its multiples $\alpha S$.
To this end
let us express the quadratic relation
by using the following notation.
For each singular fiber $F_j$,
denote by $f_j(S)$ the fraction $k_j(S)/m_j$.
Let $P(x) = x(1-x)$ for real numbers $x$.
Then the quadratic relation can be expressed as
\begin{equation}
\label{qrforS}
\sum_{j=1}^s P(f_j(S))\, m_j = 2\chi({\cal O}_X).
\end{equation}
Note again that this formula is independent of the choice of orientation
of the fibers, since changing the orientation of the components of $F_j$
simply switches the two factors $f_j$ and $1-f_j$
in each $P$ term of the above sum.

Now note that if one fixes orientations for all of the singular fibers,
and if a torsion section $S$ has component numbers $\{k_j(S)\}$
using these orientations,
then for any $\alpha$,
the component numbers for the multiple $\alpha S$ of $S$
in these orientations will be
\[
k_j(\alpha S) \equiv \alpha k_j(S) \mod m_j.
\]

Let $\langle x \rangle$ denote the smallest non-negative real number
in the residue class modulo $\Bbb Z$ of a real number $x$;
we have $\langle x \rangle \in [0,1)$ always.
The above formula can be expressed using the fractions $f_j$ by
\[
f_j(\alpha S) = \langle \alpha f_j(S) \rangle.
\]
Hence the quadratic relation (\ref{qrforS})
applied to the section $\alpha S$ is
\begin{equation}
\label{qrforalphaS}
\sum_{j=1}^s
P( \langle \alpha f_j(S) \rangle ) m_j = 2\chi({\cal O}_X).
\end{equation}

We want to reorganize this sum
to collect terms with the same $f_j$.
Let $S$ be any torsion section whose order divides $n$.
Then $nS$ is the zero-section $S_0$,
all of whose component numbers are zero;
hence for each $j$,
$n k_j(S)$ is divisible by $m_j$.
In terms of the fractions $f_j$, this is equivalent to having
$f_j(S)$ being a multiple of $1/n$ for every $j$.

Define rational numbers ${M'}_{i,n}(S)$
for $0 \leq i \leq n-1$ by the formula
\[
{M'}_{i,n}(S) =
(\sum\begin{Sb} j\text{ with }\\ f_j(S) = i /n \end{Sb} m_j)/12\chi.
\]

Roughly speaking,
${M'}_{i,n}(S)$ is the fraction of the total sum $\sum_j m_j = 12\chi$
contributed by fibers where $S$ meets the component
which is ``distance $i$'' from the component meeting the zero-section $S_0$.
(This distance is measured in units of $m_j/n$,
and is in the direction of the orientation of the fiber.)
We will call these fractions ${M'}_{i,n}(S)$
the {\em oriented distribution numbers} for the section $S$.
Obviously
\[
\sum_{i=0}^{n-1} {M'}_{i,n}(S) = 1.
\]

Reorganizing the quadratic relation (\ref{qrforalphaS}) for $\alpha S$
to sum over the $i$'s gives
\[
\sum_{i=0}^{n-1} \sum\begin{Sb} j\text{ with }\\ f_j(S) = i /n \end{Sb}
P( \langle \alpha i / n \rangle ) m_j = 2\chi({\cal O}_X),
\]
which, after dividing through by $12\chi({\cal O}_X)$,
we rewrite using the oriented distribution numbers as
\[
\sum_{i=0}^{n-1} P( \langle \alpha i / n \rangle ) {M'}_{i,n}(S)
= 1/6.
\]

If we disregard orientation,
we are led to defining rational numbers $M_{i,n}(S)$
for $0 \leq i \leq [n/2]$ by
\[
M_{i,n}(S) = \begin{cases}
{M'}_{0,n}(S) & \text{ if } i=0,\\
{M'}_{i,n}(S) + {M'}_{n-i,n}(S) & \text{ if } 1 \leq i < n/2 ,\text{ and } \\
{M'}_{n/2,n}(S) & \text{ if } i = n/2.
\end{cases}
\]
Note that
$M_{i,n}(S)$ is the fraction of the total sum $\sum_j m_j = 12\chi$
contributed by fibers where $S$ meets the component
which is ``distance $i$'' from the component meeting the zero-section $S_0$,
in either of the two orientations.
We therefore refer to the $M_{i,n}(S)$ as the
{\em unoriented distribution numbers} for the section $S$.
As with the oriented distribution numbers, we have
\begin{equation}
\label{sumMin1}
\sum_{i=0}^{[n/2]} M_{i,n}(S) = 1.
\end{equation}

Since $\langle -x \rangle = 1 - \langle x \rangle$,
and $\langle \ell + x \rangle = \langle x \rangle$ for integers $\ell$,
we see that
$\langle \alpha (n-i)/n \rangle = 1 - \langle \alpha i/n \rangle$.
Therefore since the polynomial $P(x)$ has the same value at $x$ as at $1-x$,
the quadratic relation written in terms of the oriented distribution numbers
has the coefficient of ${M'}_{i,n}(S)$ equal to the coefficient
of ${M'}_{n-i,n}(S)$.
Hence we may combine these terms, obtaining the following
relation for the unoriented distribution numbers.
(Note that the $i=0$ term has been dropped, since it is zero.)

\begin{lemma}
\label{qrforM}
If $S$ is a torsion section of order exactly $n$,
and $\alpha$ is not divisible by $n$, then
\[
\sum_{i=1}^{[n/2]} P( \langle \alpha i / n \rangle ) M_{i,n}(S)
= 1/6.
\]
\end{lemma}

Note that one gets the same equation relating the
unoriented distribution numbers $M_{i,n}(S)$
using either $\alpha$ or $n-\alpha$.
Hence the above set of equations are, a priori,
only $[n/2]$ different equations.

One can express the above set of equations in matrix form.
Let $\P_n$ be the square matrix of size $[n/2]$
whose ${\alpha i}^{th}$ entry is $P( \langle \alpha i / n \rangle )$.
Let $\1_n$ be a vector of length $[n/2]$,
which has every coordinate equal to $1$.
Finally let $\M_n(S)$ be the column vector of $M_{i,n}(S)$'s,
for $1 \leq i \leq [n/2]$.
The equations of \ref{qrforM} can then be expressed in matrix form as follows.

\begin{lemma}
\label{PMnequation}
If $S$ is a nonzero torsion section of order exactly $n$, then
\[
\P_n \M_n(S) = (1/6) \1_n.
\]
\end{lemma}

For example, suppose that $n=2$,
so that we are discussing a $2$-torsion section $S$ of $X$.
The linear system above reduces to the single equation
$(1/2) M_{1,2}(S) = 1/6$,
so that $M_{1,2}(S) = 1/3$
and therefore $M_{0,2}(S) = 2/3$ by (\ref{sumMin1}).

In case $n=3$, the linear system again reduces to the single equation
$(2/9) M_{1,3}(S) = 1/6$,
so that $M_{1,3}(S) = 3/4$
and therefore $M_{0,3}(S) = 1/4$ by (\ref{sumMin1}).

Note that these distribution numbers are independent of $S$!
This is our result in general for $p$-torsion sections,
which will be discussed in the next section.

As a last example, suppose that $n=5$.
Then the linear system above is the two equations
$(4/25)M_{1,5}(S) + (6/25)M_{2,5} = 1/6$
and
$(6/25)M_{1,5}(S) + (4/25)M_{2,5} = 1/6$
which leads to $M_{1,5}(S) = M_{2,5}(S) = 5/12$, and $M_{0,5}(S)=1/6$.
Note again that the distribution numbers are independent of $S$,
as we have seen in the previous examples.
However here we see more: the distribution numbers $M_{i,5}$
for $i\neq 0$ are {\em equal}.  This is the ``equidistribution''
property mentioned in the Introduction, and is our main result
for general primes $p$.

\section{Equidistribution for Torsion Sections of Prime Order}

In this section we will use the linear system describing the
unoriented distribution numbers $M_{i,n}(S)$
to compute these numbers for a torsion section $S$ of prime order $p$.

Our first task is to show that the matrix $\P_p$ is invertible,
which therefore shows that the distribution numbers are determined by
the linear system of Lemma \ref{PMnequation}.

\begin{proposition}
\label{P_pinvertible}
Fix an odd prime $p$.  Let $\P_p$ be the square matrix of size $(p-1)/2$
whose $\alpha\,i$ entry is $P(\langle \alpha i / p \rangle )$.
Then $\P_p$ is invertible.
\end{proposition}

\begin{pf}
Let $V$ be ${\Bbb C}^p$,
with coordinates indexed by residue classes modulo $p$.
A vector $v \in V$ will have as its $k^{th}$ coordinate the number $v_k$.
The group of units $G = {({\Bbb Z}/p)}^\times$
acts on $V$, by the formula
\[
{(a \cdot v)}_k = v_{ak \mod p}
\]
for $a \in G$.

Define vectors $S_\alpha \in V$ by setting
${(S_\alpha)}_k = P(\langle \alpha k /p \rangle)$.
We obviously have $S_\alpha = S_{p-\alpha}$.
Also note that ${(S_\alpha)}_0 = 0$ for every $\alpha$,
and ${(S_\alpha)}_{p-k} = {(S_\alpha)}_k$ for every $\alpha$ and $k$.

Therefore $S_\alpha$ lies in the subspace $V^+$ defined by
\[
V^+ = \{v \in V \,|\, v_0=0 \text{ and } v_{p-k} = v_k \text{ for all }k \}.
\]
The dimension of $V^+$ is $(p-1)/2$.
Moreover $V^+$ is stable under the $G$-action,
so that we can consider $G$ to be acting on $V^+$ also, with the same formula.
Further note that $a\cdot S_1 = S_a$ for every $a \in G$.

Note that the $\alpha^{th}$ row of $\P_p$ consists of the
coordinates ${(S_\alpha)}_k$ for $k=1,\dots,(p-1)/2$;
hence to prove that $\P_p$ is invertible,
it suffices to show that the vectors
$S_\alpha$, for $\alpha = 1,\dots,(p-1)/2$,
are independent.
Equivalently, we can show that the vectors $\{S_\alpha\}$ span $V^+$.

Let $\chi:G \to {\Bbb C}^\times$, be a character of $G$.
For $v \in V^+$, define
\[
\omega_\chi(v) = \sum_{a \in G} \chi(a) \, a \cdot v.
\]
Note that $\omega_\chi(v)$ is an eigenvector of the action of $G$ on $V^+$,
with eigenvalue $\chi^{-1}$, that is,
$a \cdot \omega_\chi(v) = \chi^{-1}(a) \omega_\chi(v)$.

Let us focus in on the special element
\[
S_\chi = \omega_\chi(S_1) = \sum_{a \in G} \chi(a) S_a.
\]
Since this vector is a linear combination of the $S_a$'s,
to show that the $S_a$'s span $V^+$,
it suffices to show that the $S_\chi$'s do.

We focus attention on the {\em even} characters $\chi$,
that is, those characters with $\chi(-1)=1$.
(If $\chi(-1)=-1$, $\chi$ is called {\em odd}.)
There are exactly $(p-1)/2$ even characters for $G$, and we will
finish the proof if we can show that
for the $(p-1)/2$ even characters $\chi$,
the vectors $S_\chi$ are non-zero;
if so, they will be linearly independent, since they
are eigenvectors for distinct characters.
Hence they will span $V^+$.

In fact we need only show that the first coordinate of each $S_\chi$
(for $\chi$ even) is non-zero.  This is the number
\[
s_\chi = {(S_\chi)}_1 = \sum_{a \in G} \chi(a) P(\langle a/p \rangle ).
\]

If $\chi$ is the trivial character which is identically $1$,
then $s_\chi$ is obviously non-zero,
so we may assume $\chi$
is one of the $(p-3)/2$ nontrivial even characters of $G$.

Let $\B_k(X)$ denote the $k^{th}$ Bernoulli polynomial.
In particular, we have
\[
\B_2(X) = X^2 - X + 1/6 = 1/6 - P(X), \text{ or }
P(X) = 1/6 - \B_2(X) = \B_2(0) - \B_2(X).
\]

Therefore the number $s_\chi$ can be written as
\[
s_\chi = \sum_{a \in G} \chi(a) [\B_2(0) - \B_2(\langle a/p \rangle ) ],
\]
and since $\chi$ is nontrivial, $\sum_{a \in G} \chi(a) = 0$,
so that
\[
s_\chi = - \sum_{a \in G} \chi(a) \B_2(\langle a/p \rangle ).
\]

Now the definition of the generalized Bernoulli numbers $B_{k,\chi}$
according to Leopoldt (see \cite[page 37]{lang}, or \cite[Proposition
4.1]{washington}), gives
\[
B_{k,\chi} = p^{k-1} \sum_{a=0}^{p-1} \chi(a) \B_k(\langle a/p \rangle);
\]
hence we have
\[
s_\chi = {-1 \over p} B_{2,\chi}
\]
for an even nontrivial character $\chi$ on $G$.

Now the classical theorem that $B_{k,\chi}$ is nonzero
when $k$ and $\chi$ have the same parity
shows that for every even character $\chi$, $s_\chi \neq 0$,
finishing the proof.
This theorem is proved by using the associated $L$-series
\[
L(s,\chi) = \sum_{n=1}^\infty {\chi(n) \over n^s}
\]
and noting first that for $n \geq 1$, $B_{n,\chi} = -n L(1-n,\chi)$
(see \cite[Proposition 16.6.2]{ireland-rosen}
or \cite[Theorem 4.2]{washington}).
Hence it suffices to show that for even $\chi$, $L(1-n,\chi)$ is nonzero.
This follows from considering the functional equation for $L$, which
for even $\chi$ can be written as
(see \cite[page 29]{washington})
\[
L(1-s,\chi) = \Gamma(s) {2 \over \tau(\chi)}{({p \over 2\pi})}^s
L(s,\overline{\chi})
\]
where $\tau(\chi) = \sum_{a=1}^p \overline{\chi}(a) e^{2\pi i a/p}$ is a Gauss
sum, and is known to be non-zero (it has absolute value $\sqrt(p)$, see
\cite[Lemma 4.8]{washington}).
Therefore to finish we must check that for even $n \geq 2$,
$L(n,\overline{\chi}) \neq 0$, and this follows since for any $\chi$,
$L(s,\chi)$ is nonzero whenever $\Re(s) >1$, by considering the
Euler product expansion for $L$, which is
\[
L(s,\chi) = \prod_p {(1 - \chi(p)p^{-s})}^{-1}.
\]
\end{pf}

The techniques of the proof given above was inspired by those used
in a similar situation in \cite[appendix to section 5]{reid},
following \cite{morrison-stevens}.
I am indebted to Miles Reid for leading me in this direction.

The first corollary of the above Proposition is that
the unoriented distribution numbers $M_{i,n}(S)$
for a torsion section $S$ of prime order
are determined by the linear system of Lemma \ref{PMnequation}.
This proves that they are independent of $S$!
Moreover, we can easily calculate them from the linear system.
We require one lemma concerning the coefficient matrix $\P_p$:

\begin{lemma}
\label{Pprowsums}
For an odd prime $p$, and for any fixed $\alpha = 1,\dots,(p-1)/2$,
\[
\sum_{i=1}^{(p-1)/2} P(\langle \alpha i /p \rangle) = { p^2-1 \over 12p}.
\]
\end{lemma}

\begin{pf}
For fixed $\alpha$, the integers $\alpha i$ as $i$ ranges from $1$ to $p-1$
themselves range over a complete set of representatives for the nonzero
residue classes modulo $p$.
Therefore, since $P(\langle \alpha i /p \rangle)$ is an even function of $i$,
the set of values of this function for $i=1,\dots, (p-1)/2$
are the same as the set of values for the function when $\alpha = 1$.
In this range we have $\langle i /p \rangle = i/p$,
so that for any $\alpha$, the sum above is simply
\[
\sum_{i=1}^{(p-1)/2} (i/p)(1-i/p) = { p^2-1 \over 12p}
\]
as claimed.
\end{pf}

Our main result is the following.

\begin{theorem}
\label{Minvalues}
Let $S$ be a torsion section of odd prime order $p$.
Then
\[
M_{i,p}(S) = 2p/(p^2-1) \; \text{ if }\; 1 \leq i \leq (p-1)/2,
\]
and
\[
M_{0,p}(S) = 1/(p+1).
\]
If $S$ is a torsion section of order $2$, then
\[
M_{1,2}(S) = 1/3 \;\text{ and }\; M_{0,2}(S) = 2/3.
\]
\end{theorem}

\begin{pf}
If $p$ is odd,
Lemma \ref{Pprowsums} shows that the matrix $\P_p$ has constant row sums.
In other words, $\P_p \1 = ((p^2-1)/ 12p)\1$,
where $\1$ is the constant vector with every coordinate $1$.
Therefore a solution to the linear system $\P_p {\bold X} = (1/6)\1$
(of which the vector $\M_p$ of $M_{i,p}(S)$'s
is a solution by Lemma \ref{PMnequation})
is given by a constant vector $(2p / (p^2-1))\1$.
Since $\P_p$ is invertible, such a solution is unique,
showing that for $i=1,\dots,(p-1)/2$, we must have
$M_{i,p}(S) = 2p/(p^2-1)$.
The value of $M_{0,p}$ is then gotten from (\ref{sumMin1}).

For $p=2$, the computations were made at the end of the last section.
\end{pf}

We will use the notation $M_{i,p}$ for $M_{i,p}(S)$ from now on.
The computation of $M_{0,p}$ was done previously in \cite{miranda-persson3}
using a quotient argument.

We view the above theorem as an ``equidistribution'' result:
for a $p$-torsion section $S$, the number of times $S$ meets
a component which is ``distance $i$'' away from the zero component $C_0$,
(measured in units of $m_j/p$),
counted with multiplicity ($m_j$ for the fiber $F_j$),
is independent of $i$, for $i \neq 0$.

As a concrete example, take the elliptic $K3$ surface
with $6$ singular fibers, having $[m_j] = [1,1,1,7,7,7]$.
Such a $K3$ surface exists, (see \cite{miranda-persson1})
and has a torsion section $S$ of order $7$ (see \cite{miranda-persson2}).
The component numbers $\{k_j(S)\}$ must be $(0,0,0,a,b,c)$,
where we can take the orientations of the $I_7$ fibers so that
$0 \leq a,b,c \leq 3$.
The equidistribution property above says that $M_{i,7} = 7/24$ for $i=1,2,3$;
hence we are forced to have $(k_j) = (0,0,0,1,2,3)$
after possibly reordering the fibers.
Since $\chi({\cal O}_X) = 2$,
the sum of the component numbers $k_j(S)$ must be $6$
by Corollary \ref{formula3}, so this is consistent.
The quadratic relation for $S$ in this case says that
$1\cdot 6/7 +2\cdot 5/7+3\cdot 4/7 = 2\chi({\cal O}_X) = 4$,
which also checks out.
Finally note that any multiple of $S$ (except $S_0$ of course)
has component numbers $\{0,0,0,1\text{ or }6,2\text{ or }5,3\text{ or }4\}$
so that all of the relations above are satisfied for all multiples.

As another example, take the Beauville surface $X_{6321}$,
which is a rational elliptic surface ($\chi({\cal O}_X)=1$)
with four singular fibers
of types $I_6$, $I_3$, $I_2$, and $I_1$ (see \cite{beauville}).
Such a surface has a finite Mordell-Weil group of order $6$
(see \cite{miranda-persson4}).
Using the relations given here for the component numbers,
we can easily deduce them.
Let $S_2$ be a torsion section of order $2$
and $S_3$ be a torsion section of order $3$.
The component numbers for $S_2$ must be
$(0\text{ or }3,0,0\text{ or }1,0)$,
and since their sum must be $4$ by Corollary \ref{formula2},
we must have $(k_j(S_2)) = (3,0,1,0)$.
Note that these numbers do not change if we change orientations,
since $S_2$ has order $2$.
The component numbers for $S_3$ must be
$(0\text{ or }2\text{ or }4,0\text{ or }1\text{ or }2,0,0)$,
and we may choose orientations so that they are
$(0\text{ or }2,0\text{ or }1,0,0)$.
By Corollary \ref{formula3}, these numbers must sum to $3$,
so that we must have $(k_j(S_3)) = (2,1,0,0)$.
In particular we have with these orientations that the sum of these
sections, which is a section $S_6$ of order $6$ generating the Mordell-Weil
group of $X$, has component numbers $(5,1,1,0)$.
(We may re-orient the fibers to achieve $(k_j(S_6)) = (1,1,1,0)$ if we desire.)

As a final example let us consider the Beauville surface $X=X_{3333}$,
which is the Hessian pencil surface obtained by taking a smooth cubic curve
$C$ in the plane and considering the pencil formed by $C$ and its Hessian $H$.
One obtains after blowing up the $9$ base points of this pencil
(which are the $9$ flexes of $C$)
a rational elliptic surface ($\chi=1$) with four $I_3$ fibers.
(One gets the same surface starting from any smooth $C$!).
The Mordell-Weil group of $X$ is of order $9$, with the zero-section
$S_0$ and $8$ other sections of order $3$.
The component numbers for a section are all $0$, $1$, or $2$,
and if $S$ is not the zero section, then
the component number sum formula Corollary \ref{formula4}
says that exactly one $k_j(S)$ is $0$, all others are non-zero.
We may choose an ordering of the fibers and orientations so that
one of the order $3$ sections $S_1$ has component numbers $(0,1,1,1)$.
Let $S_2$ be another order $3$ section, not equal to $2S_1$,
so that $S_1$ and $S_2$ generate the Mordell-Weil group.
By reordering the fibers we may assume that the component numbers for $S_2$
are $(a,0,b,c)$ where $a,b,c$ are not $0$;
after reorienting the first fiber we may further assume that $a=1$.
Since the sum $S_3=S_1+S_2$ will then have component numbers
$(1,1,b+1,c+1)$ and must have three non-zero component numbers,
we may (after possibly switching the last two fibers)
assume that $b=2$ and $c=1$,
so that $(k_j(S_2)) = (1,0,2,1)$.
This then determines the component numbers for all $8$ nonzero sections of $X$.

As an application of these computations of the unoriented distribution numbers
$M_{i,p}$ for prime $p$, we obtain a divisibility result.
Since each of the numbers $12\chi M_{i,p}(S)$ is an integer,
and since $p$ and $p^2-1$ are always relatively prime,
we obtain the following.

\begin{corollary}
\label{pdiv}
Let $p$ be an odd prime and suppose that
$f$ admits a section of order exactly $p$.
Then $(p^2-1)/2$ divides $e = 12\chi$.
\end{corollary}

This improves a divisibility result
stated in \cite[page 255]{miranda-persson3}
that $p+1$ divides $e=12\chi$, which was obtained by using $M_{0,p}$ only.
Moreover this is sharp; for the elliptic modular surface
(associated to the group $\Gamma_1(p)$)
which has a section of order $p$,
the Euler number is exactly $e = (p^2-1)/2$
(see \cite{shimura} or \cite{cox-parry}).
In some sense we can view this result as enabling us to ``divine''
the Euler number of such a universal surface, without constructing it:
if it existed, it surely would have minimal $e$.

We also have recovered what Persson and I referred to as the
``Fixed Point Rule'' in \cite[Corollary 5.5]{miranda-persson1}
(see also \cite[Lemma 3(c)]{miranda-persson2} and
\cite[Corollaries X.3.3 and X.3.4]{miranda}).
Again these works concentrated on the K3 surface case,
which has $\chi({\cal O}_X) = 2$;
the statement in this case was that
if $S$ is a $p$-torsion section on a K3 elliptic surface,
with component numbers $\{k_j\}$, then
\[
\sum \{m_j \; : \; k_j \neq 0\} = 24p/(p+1).
\]

In general of course we have (using the results of this article) that
\begin{eqnarray*}
\sum \{m_j \; : \; k_j \neq 0\} &=& 12\chi - \sum \{m_j \; : \; k_j = 0\} \\
&=& 12\chi - 12\chi\cdot M_{0,p}(S) \\
&=& 12\chi(1 - 1/(p+1)) \\
&=& 12\chi p / p+1,
\end{eqnarray*}
which generalizes the Fixed Point Rule as stated in the above articles.

Finally, I would like to mention a ``dual'' equidistribution property
which I have not been able to prove,
but which should be true.
It concerns not the indices of the components other than the
zero component which a $p$-torsion section $S$ hits,
but rather the position in $C_0^{(j)}$ which $S$ hits,
when it does hit this zero component.
Specifically, each $C_0^{(j)}$ can be identified with
${\Bbb C}^*$ by sending the nodes to $0$ and $\infty$,
and sending the point $S_0 \cap C_0^{(j)}$ to $1$.
Any torsion section $S$ of order $p$
will then hit $C_0^{(j)}$ in a $p^{th}$ root of unity.
The equidistribution property in this context should be
that each $p^{th}$ root of unity occurs equally often,
(counted properly, of course).

I do not know how to detect this root of unity
using the elementary types of techniques used in this paper.
Of course this property, properly formulated,
would be invariant under base change,
and so it would suffice to check it for the modular surfaces.
Peter Stiller has informed me that it is true
for $p \equiv 3 \mod 4$, by checking the modular surface.

\end{document}